\documentstyle[12pt]{article}
\def\stern{\mbox{\LARGE $^*$}}
\begin{document}
\includegraphics{zibheader.eps}
\vspace*{6.5cm}
\begin{center}
{\Large Wolfram Koepf}
\vspace*{3mm}\\
{\Large \hspace*{3mm}Dieter Schmersau\stern}
\vspace*{2cm}\\
\Large{{\bf Weinstein's Functions and the \\[4mm]
Askey-Gasper Identity}}
\end{center}
\vfill
\vfill
\vfill
$*$ Fachbereich Mathematik und Informatik der Freien Universit\"at Berlin
\vfill
\hrule
\vspace*{3mm}
Preprint SC 96--6 (Februar 1996)

\thispagestyle{empty}
\setcounter{page}{0}
\eject
\begin{tabular}{l}
\\[20cm]
\end{tabular}
\thispagestyle{empty}
\setcounter{page}{0}
\hoffset -1cm
\footskip=1cm
\parindent=0pt
\title{Weinstein's Functions and the \\[1mm]
Askey-Gasper Identity}
\date{}
\author{Wolfram Koepf\\
Dieter Schmersau\\
email: {\tt koepf@zib-berlin.de}
}
\maketitle
%
%
%
%
\setlength{\textheight}{8.5in}
\setlength{\textwidth}{6.5in}
\setlength{\evensidemargin} {0in}
\setlength{\oddsidemargin} {0in}
\setlength{\topmargin} {0in}

\def\sumprime{\mathop{\sum\nolimits'}\limits}
\newcommand{\ded}[1]{\frac{\de}{\de #1}}
\newcommand{\dedn}[2]{\frac{\de^{#2}}{\de #1^{#2}}}
\newcommand{\hypergeom}[5]{\mbox{$
_#1 F_#2\left.
\!\!
\left(
\!\!\!\!
\begin{array}{c}
\multicolumn{1}{c}{\begin{array}{c}
#3
\end{array}}\\[1mm]
\multicolumn{1}{c}{\begin{array}{c}
#4
            \end{array}}\end{array}
\!\!\!\!
\right| \displaystyle{#5}\right)
$}
}

\newcommand{\punktpunkt}[1]{\stackrel{\mbox{..}}{#1}\!}

\newcommand{\RR}{{\rm I\! R}}
\newcommand{\DD}{{\rm I\! D}} 
\newcommand{\CC}{{\; \rm l\!\!\!	C}}
\newcommand{\NN}{{\rm I\!I\!\! N}}
\newcommand{\ZZ}{\Bbb Z}
\newcommand{\C}{{\rm {\mbox{C{\llap{{\vrule height1.5ex}\kern.4em}}}}}}
\newcommand{\N} {{\rm {\mbox{\protect\makebox[.15em][l]{I}N}}}}
\renewcommand{\H} {{\rm {\mbox{\protect\makebox[.15em][l]{I}H}}}}
\newcommand{\Q} {{\rm {\mbox{Q{\llap{{\vrule height1.5ex}\kern.5em}}}}}}
\newcommand{\R} {{\rm {\mbox{\protect\makebox[.15em][l]{I}R}}}}
\newcommand{\D} {{\rm {\mbox{\protect\makebox[.15em][l]{I}D}}}}
\newcommand{\Z} {{\rm {\mbox{\protect\makebox[.2em][l]{\sf Z}\sf Z}}}}
\newcommand{\Rp}
{{\rm {\mbox{$\mbox{\protect\makebox[.15em][l]{I}R}^{\scriptscriptstyle+}$}\index{R+@\mbox{$\mbox{\protect\makebox[.15em][l]{I}R}^{\scriptscriptstyle+}$}, 
Notation for the positive real numbers}}}}
\newcommand{\Rm}
{{\rm {\mbox{$\mbox{\protect\makebox[.15em][l]{I}R}^{\scriptscriptstyle-}$}\index{R-@\mbox{$\mbox{\protect\makebox[.15em][l]{I}R}^{\scriptscriptstyle-}$}, 
Notation for the negative real numbers}}}}
\newcommand{\Pt}{\widetilde{P}}
\newcommand{\Bt}{\widetilde{B}}
\newcommand{\Cd}{\widehat{{\rm \;l\!\!\! C}}}
\newcommand{\bq}[1]{|#1|^2}
\newcommand{\vf}[1]{(1-\bq{#1})}
\newcommand{\vfq}[1]{\vf{#1}^2}
\newcommand{\Bf}[1]{\vf{#1}\left|\frac{f''}{f'}(#1)\right|}
\newcommand{\Nf}[1]{\vfq{#1}|S_f(#1)|}
\newcommand{\Kf}[1]{\left|-\kon{#1}+\frac{1}{2}\vf{#1}\frac{f''}{f'}(#1)\right|}
\newcommand{\ed}[1]{\frac{1}{#1}}
\newcommand{\aut}[1]{\frac{{\textstyle{z+#1}}}{{\textstyle{1+\kon{#1}z}}}}
\newcommand{\au}[2]{#1\aut{#2}}
\newcommand{\subs}[2]{\left. \makebox{\rule{0in}{2.5ex}} #1 \rb_{#2}}
\newcommand{\subst}[3]{\left. \makebox{\rule{0in}{2.5ex}} #1 \rb_{#2}^{#3}}
\newcommand{\funkdef}[3]{\left\{\begin{array}{ccc}
                                #1 && \mbox{\rm{if} $#2$ } \\
                                #3 && \mbox{\rm{otherwise}}
                                \end{array}  
                         \right.}       
\newcommand{\funkdeff}[4]{\left\{\begin{array}{ccc}
                                 #1 && \mbox{\rm{if} $#2$ } \\
                                 #3 && \mbox{\rm{if} $#4$ } 
                                 \end{array}
                          \right.}
\newcommand{\funkdefff}[6]{\left\{\begin{array}{ccc}
                                 #1 && \mbox{{if} $#2$ } \\
                                 #3 && \mbox{{if} $#4$ } \\
                                 #5 && \mbox{{if} $#6$ }
                                 \end{array}
                          \right.}
\newcommand{\ueber}[2]{
                       \Big( \!
                       {{\small
                       \begin{array}{c}
                          #1\\
                          #2
                          \end{array}
                       }}
                       \! \Big) }
\newcommand{\function}[4]{
                          \begin{array}{rcl}#1&\pf&#2\\
                          #3&\mapsto &#4
                          \end{array} }
\newcommand{\pr}{\vspace{-2mm}\absatz{{\sl Proof:}}\hspace{5mm}}
\newcommand{\eop}{\hfill$\Box$\par\medskip\noindent}
\newcommand{\absatz}{\par\medskip\noindent}
\renewcommand{\Re}{{\rm Re\:}}
\renewcommand{\Im}{{\rm Im\:}}
\newcommand{\co}{{\rm co\:}}
\newcommand{\coq}{\overline{{\rm co}} \:}
\newcommand{\ex}{{\rm E\:}}
\newcommand{\In}{\in}
\newcommand{\ro}{\varrho}
\newcommand{\om}{\omega}
\newcommand{\al}{\alpha}
\newcommand{\bb}{\beta}
\newcommand{\la}{\lambda}
\newcommand{\eps}{\varepsilon}
\newcommand{\ph}{\varphi}
\renewcommand{\phi}{\varphi}
\newcommand{\si}{\sigma}
\newcommand{\ka}{\varkappa}
\newcommand{\th}{\theta}
\newcommand{\g}{\gamma}
\newcommand{\de}{\partial}
\newcommand{\fD}{f(\D)}
\newcommand{\sumi}{\sum\limits_{k=0}^{\infty}}
\newcommand{\sumei}{\sum\limits_{k=1}^{\infty}}
\newcommand{\union}{\bigcup\limits_{k=1}^{n}}
\newcommand{\sumn}{\sum\limits_{k=1}^{n}}
\newcommand{\prodn}{\prod\limits_{k=1}^{n}}
\newcommand{\intd}{\int\limits_{\de\DD}}
\newcommand{\menge}[3]{\left\{#1 \In #2 \; \lb \; #3 \right. \right\} }
\newcommand{\mk}{\mu_{k}}
\newcommand{\xk}{x_k}
\newcommand{\yk}{y_k}
\newcommand{\xn}{x_n}
\newcommand{\yn}{y_n}
\newcommand{\ak}{\al_k}
\newcommand{\bk}{\bb_k}
\newcommand{\kn}{\mbox{$(k=1, \ldots ,n)$}}
\newcommand{\kno}{\mbox{$k=1, \ldots ,n$}}
\newcommand{\sub}{\prec}
\newcommand{\ld}[1]{\frac{f''}{f'}(#1)} 
\newcommand{\limn}{\lim\limits_{n\rightarrow\infty}}
\newcommand{\lsz}{\limsup\limits_{z\rightarrow\de\DD}}
\newcommand{\limr}{\lim\limits_{r\rightarrow 1}}
\newcommand{\liz}{\liminf\limits_{z\rightarrow\de\DD}}
\newcommand{\supD}[1]{\sup\limits_{#1\In \DD}}
\newcommand{\infD}[1]{\inf\limits_{#1\In \DD}}
\newcommand{\maxn}{\max\limits_{1 \leqq k \leqq n}}
\newcommand{\minn}{\min\limits_{1 \leqq k \leqq n}}
\newcommand{\ord}{{\rm ord\:}}
\newcommand{\gleich}[1]{\stackrel{{{\rm (#1)}}}{\longeq}}
\newcommand{\folgt}[1]{\stackrel{{{\rm (#1)}}}{\Pf}}
\newcommand{\gegen}[1]{\stackrel{{{\rm (#1)}}}{\pf}}
\newcommand{\nach}[2]{(#1)$\dpf$(#2):}
\newcommand{\Exp}{\subset\!\subset}
\newcommand{\lleq}{\stackrel{_{{\scriptscriptstyle \Exp}}}
           {_{{\scriptscriptstyle \sim}}}}
\newcommand{\Sub}{{\rm Sub\:}}
\newcommand{\norm}[2]{\frac{#1\circ #2-#1\circ #2(0)}{(#1\circ #2)'(0)}} 
\renewcommand{\dim}{{\rm{dim}}_{_{H^{p}}}}
\renewcommand{\ll}{<\!<}
\newcommand{\1}{{\bf{1}}}
\newcommand{\2}{{\bf{2}}}
\newcommand{\3}{{\bf{3}}}
\newcommand{\4}{{\bf{4}}}
\newcommand{\5}{{\bf{5}}}
\newcommand{\6}{{\bf{6}}}
\newcommand{\7}{{\bf{7}}}
\newcommand{\8}{{\bf{8}}}
\newcommand{\9}{{\bf{9}}}
\newcommand{\0}{{\bf{0}}}
\def\finis{\hbox{$\bigtriangleup$}}    


\newcommand{\luc}{locally uniform convergence}
\newcommand{\hp}{half\-plane}
\newcommand{\sq}{sequence}
\newcommand{\an}{analytic}
\newcommand{\af}{analytic\ function}
\newcommand{\Ne}{Ne\-ha\-ri\ ex\-pression}
\newcommand{\Ke}{Koe\-be\ ex\-pression}
\newcommand{\Be}{Becker\ ex\-pression}
\newcommand{\fc}{function}
\newcommand{\uv}{univalent}
\newcommand{\uf}{univalent\ function}
\newcommand{\SC}{Schwarz-Chri\-stof\-fel}
\newcommand{\pol}{polyg\-on}
\newcommand{\cv}{con\-vex}
\newcommand{\ctc}{close-to-con\-vex}
\newcommand{\Ck}{Ca\-ra\-th\'eo\-dory\ kernel}
\newcommand{\st}{sector}
\newcommand{\sm}{similar}
\newcommand{\br}{boundary\ rotation}
\newcommand{\bbr}{bounded\ \br}
\newcommand{\KM}{Krein-Mil\-man}


\newcommand{\til}{\widetilde}
\newcommand{\pf}{\rightarrow}
\newcommand{\Pf}{\;\;\;\longrightarrow\;\;\;}
\newcommand{\dpf}{\Rightarrow}
\newcommand{\Dpf}{\;\;\;\Longrightarrow\;\;\;}
\newcommand{\kon}{\overline}
\newcommand{\be}{\begin{equation}}
\newcommand{\ee}{\end{equation}}
\newcommand{\bea}{\begin{eqnarray}}
\newcommand{\eea}{\end{eqnarray}}
\newcommand{\beao}{\begin{eqnarray*}}
\newcommand{\eeao}{\end{eqnarray*}}
\newcommand{\lequiv}{\;\;\;\Longleftrightarrow\;\;\;}
\newcommand{\longeq}{\;\;\;=\!\!=\!\!=\;\;\;}
\newcommand{\leqq}{\leq}
\newcommand{\geqq}{\geq}
\newcommand{\gl}{\;\leftrightarrow\;}
\newcommand{\lk}{\left(}
\newcommand{\rk}{\right)}
\newcommand{\lb}{\left|}
\newcommand{\rb}{\right|}
\newcommand{\bi}{\bibitem}


\newcommand{\bT}{\begin{theorem}}
\newcommand{\eT}{\end{theorem}}
\newcommand{\bL}{\begin{lemma}}
\newcommand{\eL}{\end{lemma}}
\newcommand{\bC}{\begin{corollary}}  
\newcommand{\eC}{\end{corollary}}
\newcommand{\bt}{\begin{tabbing} 12345 \= \kill}
\newcommand{\et}{\end{tabbing}}


\hyphenation{qua-si-disk qua-si-circle 
non-smooth
pa-ram-e-trized pa-ram-e-tri-zation
geo-met-ric
}


\newcommand{\bbegin}{{\bf{begin}}}
\newcommand{\eend}{{\bf{end}}}
\newcommand{\iif}{{\bf{if}}}
\newcommand{\tthen}{{\bf{then}}}
\newcommand{\wwhile}{{\bf{while}}}
\newcommand{\ddo}{{\bf{do}}}
\newcommand{\ffor}{{\bf{for}}}
\newcommand{\sstep}{{\bf{step}}}
\newcommand{\llet}{{\bf{let}}}
\newcommand{\pprocedure}{{\bf{procedure}}}
\newcommand{\aand}{{\bf{and}}}
\newcommand{\nnot}{{\bf{not}}}
\newcommand{\oor}{{\bf{or}}}
\newcommand{\lllet}{{\bf{let}}}
\newcommand{\eexit}{{\bf{exit}}}
\newcommand{\rreturn}{{\bf{return}}}
\newcommand{\uuntil}{{\bf{until}}}

\newcommand{\abs}{\\[3mm]}

\newtheorem{theorem}{Theorem}
\newtheorem{algorithm}{Algorithm}
\newtheorem{lemma}{Lemma}
\newtheorem{corollary}[theorem]{Corollary}
\newtheorem{definition}{Definition}
\begin{abstract}
In his 1984 proof of the Bieberbach and Milin conjectures de Branges used
a positivity result of special functions which follows from an
identity about Jacobi polynomial sums that was found by Askey and Gasper in
1973, published in 1976.

In 1991 Weinstein presented another proof of the Bieberbach and Milin 
conjectures, also using a special function system which (by Todorov and Wilf) 
was realized to be the same as de Branges'.

In this article, we show how a variant of the Askey-Gasper identity can be
deduced by a straightforward examination of Weinstein's functions
which intimately are related with a L\"owner chain of the Koebe function, 
and therefore with univalent functions.
\end{abstract}
 
\section{Introduction}

Let $S$ denote the family of analytic and univalent functions
$f(z)=z+a_2 z^2+\ldots$ of the unit disk $\D$. $S$ is compact with
respect to the topology of locally uniform convergence so that
\mbox{$k_{n}:=\max\limits_{{f\in S}}|a_{n}(f)|$} exists. In 1916
Bieberbach \cite{Bieberbach}
proved that $k_2=2$, with equality if and only if $f$ is a rotation of
the {\sl Koebe function}
\be
K(z):=\frac{z}{(1-z)^{2}}=\ed{4}\lk\lk\frac{1+z}{1-z}\rk^{2}-1\rk=
\sum\limits_{n=1}^{\infty}nz^{n}\;,
\label{eq:Koebe function}
\ee
and in a footnote he mentioned ``Vielleicht ist
\"uberhaupt $k_{n}=n$.''. This statement is known as the {\sl Bieberbach
conjecture.}

In 1923 L\"owner \cite{Loewner2} proved the Bieberbach conjecture for
$n=3$. His method was to embed a univalent function $f(z)$ into a {\sl L\"owner
chain\/}, i.e.\ a family $\left\{f(z,t) \; \lb \; t\geqq 0 \right.
\right\}$ of univalent functions of the form
\[
f(z,t)=e^{t}z+\sum\limits_{n=2}^{\infty}a_{n}(t)z^{n},
\;\;\;(z\in\D, t\geqq 0, a_{n}(t)\in\C\;(n\geqq 2))
\]
which start with $f$
\[
f(z,0)=f(z)\;,
\]
and for which the relation
\be
\Re p(z,t)=\Re\lk\frac{{\dot{f}}(z,t)}{z f'(z,t)}\rk>0
\;\;\;\;\;\;\;\;\;(z\in\D)
\label{eq:p(z,t)}
\ee
is satisfied. Here $'$ and $\dot{}\:$ denote the partial derivatives
with respect to $z$ and  $t$, respectively. Equation~(\ref{eq:p(z,t)})
is referred to as the
{\sl L\"owner differential equation\/}, and
geometrically it states that the image domains of
$f_t$ expand as $t$ increases.

The history of the Bieberbach conjecture showed that it was easier to
obtain results about the {\sl logarithmic coefficients\/} of a univalent
function $f$, i.e.\ the coefficients $d_n$ of the expansion
\[
\phi(z)=\ln\frac{f(z)}{z}=:\sum\limits_{n=1}^{\infty}d_{n}z^{n}
\]
rather than for the coefficients $a_n$ of $f$ itself. So Lebedev and
Milin \cite{LM}
in the mid sixties developed methods to exponentiate such
information. They proved that if
for $f\in S$ the {\sl Milin conjecture\/}
\[
\sum\limits_{k=1}^{n}(n+1-k)\lk k|d_{k}|^{2}-\frac{4}{k}\rk\leqq 0
\]
on its logarithmic coefficients is satisfied for some $n\in\N$, then 
the Bieberbach conjecture for the index $n+1$ follows.

In 1984 de Branges \cite{Bra}
verified the Milin, and therefore the Bieberbach
conjecture, and in 1991, Weinstein \cite{Weinstein} gave a different proof.
A reference other than \cite{Bra} concerning de Branges' proof is \cite{FP},
and a German language summary of the history of the Bieberbach conjecture
and its proofs was given in \cite{Koepf}.

Both proofs use the positivity of special function systems, and
independently Todorov \cite{Todorov} and Wilf \cite{Wilf}
showed that both de Branges' and Weinstein's functions essentially 
are the same (see also \cite{KSdeBranges}),
\be
\dot{\tau_{k}^n}(t)=-k\Lambda_k^n(t)
\;,
\label{eq:deBrangesWeinstein}
\ee
$\tau_k^n(t)$ denoting the de Branges functions and $\Lambda_k^n(t)$
denoting the Weinstein functions, respectively.

Whereas de Branges applied an identity of Askey and Gasper \cite{AG}
to his function system, Weinstein applied an addition theorem for
Legendre polynomials to his function system to deduce the positivity
result needed.

The identity of Askey and Gasper used by de Branges
was stated in (\cite{AG}, (1.16)) in the form
\be
\sum_{j=0}^{n} P_j^{(\alpha,0)}(x)=
\sum_{j=0}^{[n/2]}
\frac{(1/2)_j\,\lk\frac{\al+2}{2}\rk_{n-j}\,\lk\frac{\al+3}{2}\rk_{n-2j}(n-2j)!}
{j!\,\lk\frac{\al+3}{2}\rk_{n-j}\,\lk\frac{\al+1}{2}\rk_{n-2j}(\al+1)_{n-2j}}
\lk C_{n-2j}^{(\al+1)/2}\lk\sqrt{\frac{1+x}{2}}\rk\rk^2
\;,
\label{eq:AG Identity}
\ee
where $C_n^\lambda(x)$ denote the Gegenbauer polynomials,
$P_j^{(\alpha,\beta)}(x)$ denote the Jacobi polynomials (see e.g.\
\cite{AS}, \S~22), and
\[
(a)_j:=a(a+1)\cdots(a+j-1)=\frac{\Gamma(a+j)}{\Gamma(a)}
\]
denotes the shifted factorial (or Pochhammer symbol).

In this article, we show how a variant of the Askey-Gasper identity can be
deduced by a straightforward examination of Weinstein's functions
which intimately are related with the bounded
L\"owner chain of the Koebe function. 

The application of an addition theorem for the Gegenbauer polynomials
quite naturally arises 
in this context. We present a simple proof of 
this result so that this article is self-contained.

\section{The L\"owner Chain of the Koebe Function and the Weinstein Functions}
\label{sec:The Lowner Chain of the Koebe Function and the Weinstein Functions}

We consider the L\"owner chain
\be
w(z,t):=K^{-1}\Big(e^{-t}K(z)\Big)\quad(z\in\D, t\geq 0)
\label{eq:Koebe}
\ee
of bounded univalent functions in the unit disk $\D$
which is defined in terms of the Koebe function (\ref{eq:Koebe function}).
Since $K$ maps the unit disk onto the entire plane slit along the
negative $x$-axis in the interval $(-\infty,1/4]$,
the image $w(\D,t)$ is the unit disk with a radial slit on the
negative $x$-axis increasing with $t$.

Weinstein \cite{Weinstein} 
used the L\"owner chain (\ref{eq:Koebe}),
and showed the validity of Milin's conjecture if for all $n\geq 2$ 
the {\sl Weinstein functions} $\Lambda_k^n:\R^+\rightarrow\R
\;(k=0,\ldots,n)$ defined by
\be
\frac{e^{t}w(z,t)^{k+1}}{1-w^2(z,t)}=:
\sum\limits_{n=k}^{\infty}\Lambda_k^n(t)z^{n+1}
=W_k(z,t)
\label{eq:W(z,t)}
\;,
\ee
satisfy the relations
\begin{equation}
\Lambda_k^n(t)\geq 0\;\;\;\;\;\;(t\in\R^+,\;\;\;0\leq k\leq n)
\;.
\label{eq:pos}
\end{equation}
Weinstein did not identify the functions $\Lambda_k^n(t)$, but was
able to prove (\ref{eq:pos}) without an explicit representation.

In this section we apply Weinstein's following interesting observation 
to show that $\Lambda_k^n(t)$ are the Fourier coefficients of a
function that is connected with the Gegenbauer
and Chebyshev polynomials.

The range of the function $w=K^{-1}(e^{-t}K)$ is the unit disk
with a slit on the negative real axis. Since for all 
$\g\In\R, \g\neq 0\;({\rm mod\:\pi})$ the mapping 
\[
h_{\g}(z):=\frac{z}{1-2\cos\g\cdot z+z^{2}}
\]
maps the unit disk 
onto the unit disk with two slits on the real axis, we can interpret 
$w$ as composition $w=h_{\th}^{-1}(e^{-t}h_{\g})$ for a suitable
pair $(\th,\g)$, and a simple calculation shows that the relation
\be
\cos\g=(1-e^{-t})+e^{-t}\cos\th
\label{eq:cosg}
\ee
is valid. 
We get therefore
\bea
h_{\g}(z)&=&e^{t}\cdot h_{\th}(w(z,t))=
\frac{e^{t}w}{1-w^{2}}
\lk\frac{1-w^{2}}{1-2\cos\th\cdot w+w^{2}}\rk
\nonumber
\\
&=&
\frac{e^{t}w}{1-w^{2}}
\lk 1 + 2\sum\limits_{k=1}^{\infty}w^{k}\cos k\th\rk
=
W_0(z,t)+2\sum_{k=1}^\infty W_k(z,t)\cos k\th
\nonumber
\\
&=&W_0(z,t)
+ 2\sum\limits_{k=1}^{\infty}\lk\sum\limits_{n=k}^{\infty}
\Lambda_k^n(t)z^{n+1}\rk\cos k\th
\label{eq:hgamma}
\;.
\eea
It is easily seen that (\ref{eq:hgamma}) remains valid for the pair
$(\theta,\g)=(0,0)$, corresponding to the representation
\[
K(z)=W_0(z,t)+2\sum_{k=1}^\infty W_k(z,t)
\;.
\]
Since on the other hand $h_\g(z)$ has the Taylor expansion
\[
h_\g(z)=
\frac{z}{1-2\cos\g\cdot z+z^{2}}=
\sum_{n=0}^\infty
\frac{\sin (n+1)\g}{\sin\g}z^{n+1}
\;,
\]
equating the coefficients of $z^{n+1}$ in (\ref{eq:hgamma}) we get the
identity
\[
\frac{\sin (n+1)\g}{\sin\g}=
\Lambda_0^n(t)+2\sum_{k=1}^n \Lambda_k^n(t) \cos k\th
\;.
\]
Hence we have discovered (see also \cite{Wilf}, (2))
\bT
[Fourier Expansion]
\label{th:Weinstein Fourier}
{\rm
The Weinstein functions $\Lambda_k^n(t)$ 
satisfy the functional equation
\be
U_n\Big((1-e^{-t})+e^{-t}\cos\th\Big)
=
C_n^{1}\Big((1-e^{-t})+e^{-t}\cos\th\Big)
=
\Lambda_0^n(t)+2\sum_{k=1}^n \Lambda_k^n(t) \cos k\th
\;,
\label{eq:Weinstein Fourier}
\ee
where $U_n(x)$ denote the Chebyshev polynomials of the second kind.
}
\eT
\pr
This is an immediate consequence of the identity
\[
C_n^{1}(\cos\g)=U_n(\cos\g)=\frac{\sin (n+1)\g}{\sin\g}
\]
(see e.g.\ \cite{AS}, (22.3.16), (22.5.34)).
\eop

\section{The Weinstein Functions as Jacobi Polynomial Sums}
\label{sec:The Weinstein Functions as Jacobi Polynomial Sums}

In this section, we show that the Weinstein functions $\Lambda_k^n(t)$
can be represented as Jacobi polynomial sums.
\bT
[Jacobi Sum]
\label{th:Weinstein Jacobisum}
{\rm
The Weinstein functions have the representation
\be
\Lambda_k^n(t)=e^{-kt}\sum_{j=0}^{n-k} P_j^{(2k,0)}(1-2e^{-t})
\;,\quad\quad
(0\leq k\leq n)
\;.
\label{eq:Weinstein Jacobisum}
\ee
}
\eT
\pr
A calculation shows that $w(z,t)$ has the explicit representation
\be
w(z,t)=\frac{4e^{-t}z}{\lk 1-z+\sqrt{1-2xz+z^2}\rk^2}
\label{eq:w(z,t) explicit}
\;.
\ee
Here we use the abbreviation $x=1-2e^{-t}$. Furthermore, from
\[
W_0(z,t)=\frac{e^tw}{1-w^2}=K(z)\,\frac{1-w}{1+w}
\;,
\]
we get the explicit representation
\be
W_0(z,t)=\frac{z}{1-z}\,\frac{1}{\sqrt{1-2xz+z^2}}
\label{eq:W0}
\ee
for $W_0(z,t)$. By the definition of $W_k(z)$, we have moreover
\[
W_{k}(z,t)=\frac{e^{t}w^{k+1}}{1-w^{2}}=w^k\,W_{0}(z,t)
\;.
\]
Hence, by (\ref{eq:w(z,t) explicit})--(\ref{eq:W0}) we deduce
the explicit representation
\be
W_{k}(z,t)=e^{-kt}\,\frac{z^{k+1}}{1-z}\,
\frac{4^k}{\sqrt{1-2xz+z^2}}\,\frac{1}{\left(1-z+\sqrt{1-2xz+z^2}\right)^{2k}}
\label{eq:Wkexplicit}
\ee
for $W_k(z,t)$.

Since the Jacobi polynomials $P_j^{(\al,\beta)}(x)$
have the generating function
\be
\sum_{j=0}^\infty P_j^{(\al,\beta)}(x)\,z^j=
\frac{2^{\al+\beta}}{\sqrt{1\!-\!2xz\!+\!z^2}}\,
\frac{1}{\left(1-z+\sqrt{1\!-\!2xz\!+\!z^2}\right)^{\alpha}}\,
\frac{1}{\left(1+z+\sqrt{1\!-\!2xz\!+\!z^2}\right)^{\beta}}
\label{eq:generating function Jacobi}
\ee
(see e.g.\ \cite{AS}, (22.9.1)), comparison with (\ref{eq:Wkexplicit}) yields
\[
W_k(z,t)=e^{-kt}\,\frac{z^{k+1}}{1-z}\,
\sum_{j=0}^\infty P_j^{(2k,0)}(x)\,z^j
\;.
\]
Using the Cauchy product
\[
\ed{1-z}\,\sum_{j=0}^\infty P_j^{(2k,0)}(x)\,z^j=
\sum_{n=0}^\infty\sum_{j=0}^n P_j^{(2k,0)}(x)\,z^n
\;,
\]
we finally have
\[
W_k(z,t)
=
e^{-kt}\,z^{k+1}\,\sum_{n=0}^\infty\sum_{j=0}^n P_j^{(2k,0)}(x)\,z^n
=
\sum_{n=k}^\infty \Lambda_k^n(t)\,z^{n+1}=
\sum_{n=0}^\infty \Lambda_k^{n+k}(t)\,z^{n+k+1}
\;.
\]
Equating coefficients gives the result.
\eop

\section{Askey-Gasper Inequality for the Weinstein Functions}
\label{sec:Askey-Gasper Inequality for the Weinstein Functions}

We would like to utilize the Fourier expansion (\ref{eq:Weinstein Fourier}) of
Theorem~\ref{th:Weinstein Fourier} to
find new representations for the Weinstein functions, hence by 
Theorem~\ref{th:Weinstein Jacobisum} for the Jacobi polynomial sum
on the left hand side of (\ref{eq:AG Identity}).
Hence, we have the need to find a representation for
$C_n^{1}\Big((1-e^{-t})+e^{-t}\cos\th\Big)$.

We do a little more, and give a representation for
\be
C_n^{1}\Big(xy+\sqrt{1-x^2}\sqrt{1-y^2}\,\zeta\Big)
\;,
\label{eq:Cn1ansatz}
\ee
from which the above expression is the special case $x=y=\sqrt{1-e^{-t}},
\zeta=\cos\th$. Actually, in the next section, an even more general 
expression is considered, 
see Theorem~\ref{th:Addition Theorem for the Gegenbauer Polynomials}.
Here we outline the deduction for our particular case.

The function given by (\ref{eq:Cn1ansatz}) as a function of the
variable $\zeta$ is a polynomial of degree $n$. Hence it can be
expanded by Gegenbauer polynomials $C_j^\lambda(\zeta)\;(j=0,\ldots,n)$.
We choose $\lambda=1/2$, i.e.\ we develop in terms of Legendre 
polynomials $P_j(\zeta)=C_j^{1/2}(\zeta)$ (see e.g.\ \cite{AS}, (22.5.36)), 
\be
C_n^{1}\Big(xy+\sqrt{1-x^2}\sqrt{1-y^2}\,\zeta\Big)=
\sum_{m=0}^n A_m^n(x,y)\,C_m^{1/2}(\zeta)
\label{eq:Cn1sum}
\ee
with $A_j^n$ depending on $x$ and $y$. By the orthogonality of the
Gegenbauer polynomials,
\[
\int\limits_{-1}^1 C_j^{1/2}(\zeta)\,C_m^{1/2}(\zeta)\,d\zeta=
\funkdef{\frac{2}{2j+1}}{j=m}{0}
\;,
\]
multiplying (\ref{eq:Cn1sum}) by $C_j^{1/2}(\zeta)$, and
integrating from $\zeta=-1$ to $\zeta=1$, we get therefore
\be
A_j^n(x,y)=\frac{2j+1}{2}\int\limits_{-1}^1 
C_n^{1}\Big(xy+\sqrt{1-x^2}\sqrt{1-y^2}\,\zeta\Big)
C_j^{1/2}(\zeta)\,d\zeta
\;.
\label{eq:Ajn}
\ee
To eliminate the second (oscillating) factor $C_j^{1/2}(\zeta)$, we
utilize the identity
\be
\int\limits_{-1}^1
f(\zeta)\,C_j^\lambda(\zeta)(1-\zeta^2)^{\lambda-1/2}\,d\zeta
=\frac{2^j}{j!}\frac{\Gamma(j+\lambda)\Gamma(j+2\lambda)}
{\Gamma(\lambda)\Gamma(2j+2\lambda)}
\int\limits_{-1}^1
f^{(j)}(\zeta)\,(1-\zeta^2)^{\lambda+j-1/2}\,d\zeta
\;,
\label{eq:identity f}
\ee
which is valid for any $j$ times continuously differentiable function $f$,
and which can easily be proved by iterative partial integration 
(see e.g.\ \cite{Hua}, Chapter VII, p.\ 140).
Choosing $\lambda=1/2$ and
\[
f(\zeta):=C_n^{1}\Big(xy\!+\!\sqrt{1\!-\!x^2}\sqrt{1\!-\!y^2}\,\zeta\Big)
\]
we get (with the Gamma duplication formula (\ref{eq:duplication formula}))
\be
\int\limits_{-1}^1\!
C_n^{1}\Big(xy\!+\!\sqrt{1\!-\!x^2}\sqrt{1\!-\!y^2}\zeta\Big)\,
C_j^{1/2}(\zeta)\,d\zeta\!=\!
\frac{1}{2^j\,j!}\!\int\limits_{-1}^1\!
(1-\zeta^2)^j\,\frac{d^j}{d\zeta^j}
C_n^{1}\Big(xy\!+\!\sqrt{1\!-\!x^2}\sqrt{1\!-\!y^2}\zeta\Big)\,d\zeta
.
\label{eq:Cn1Big}
\ee
Since furthermore
\be
\frac{d^j}{d\zeta^j} C_n^\nu(\zeta)=2^j\,(\nu)_j\,C_{n-j}^{\nu+j}(\zeta)
\label{eq:derivativeidentity}
\ee
(see e.g.\ \cite{Tri}, p. 179), we get moreover
\be
\frac{1}{2^j\,j!}\int\limits_{-1}^1
(1-\zeta^2)^j\,\frac{d^j}{d\zeta^j}
C_n^{1}\Big(xy+\sqrt{1-x^2}\sqrt{1-y^2}\,\zeta\Big)\,d\zeta
=
(1-x^2)^{j/2}\,(1-y^2)^{j/2}\,Q_j^n(x,y)
\label{eq:third step}
\ee
with
\[
Q_j^n(x,y):=
\int\limits_{-1}^1
(1-\zeta^2)^j\,
C_{n-j}^{j+1}\Big(xy+\sqrt{1-x^2}\sqrt{1-y^2}\,\zeta\Big)
\,d\zeta
\;.
\]
Now observe that $Q_j^n(x,y)$ is a polynomial in the variables $x$ and $y$,
of degree $n-j$ each. In the next section we will show that the 
integral $Q_j^n(x,y)$ has zeros at both the zeros
of $C_{n-j}^{j+1}(x)$ and $C_{n-j}^{j+1}(y)$, hence, as a polynomial
of degree $n-j$ in $x$ and $y$ respectively, must be a multiple of the product
$C_{n-j}^{j+1}(x)\,C_{n-j}^{j+1}(y)$. An initial value gives
\be
Q_j^n(x,y)=
\frac{2^{2(j+1)}\,j!^2\,(n-j)!}{2(n+j+1)!}
C_{n-j}^{j+1}(x)\,C_{n-j}^{j+1}(y)
\;.
\label{eq:fourth step}
\ee
Note that the complete proof of a generalization of statement 
(\ref{eq:Cn1sum})/(\ref{eq:fourth step}) will be given in the next section.

Therefore finally, combining (\ref{eq:Ajn})--(\ref{eq:fourth step}),
we have discovered the identity
\be
A_j^n(x,y)=(2j+1)
\frac{2^{2j}\,j!^2\,(n-j)!}{(n+j+1)!}\,
(1-x^2)^{j/2}\,(1-y^2)^{j/2}\,C_{n-j}^{j+1}(x)\,C_{n-j}^{j+1}(y)
\;.
\label{eq:final Ajn}
\ee
As a first step this leads to the following Askey-Gasper type representation
for the Fourier series (\ref{eq:Weinstein Fourier}).
\bT
\label{th:Askey-Gasper type}
{\rm
The Fourier series (\ref{eq:Weinstein Fourier}) has the representation
\bea
C_n^{1}\Big((1-e^{-t})+e^{-t}\cos\th\Big)
&=&
\sum_{j=0}^n A_j^n\left(\sqrt{1-e^{-t}},\sqrt{1-e^{-t}}\right)\,
C_j^{1/2}(\cos\th)
\label{eq:Fourier by Legendre}
\\&=&
\sum_{j=0}^n 
(2j+1)
\frac{4^{j}\,j!^2\,(n-j)!}{(n+j+1)!}\,
e^{-jt}\,\lk C_{n-j}^{j+1}\left(\sqrt{1-e^{-t}}\right)\rk^2\,
P_j(\cos\th)
\;.
\nonumber
\eea
}
\eT
\pr
Set $x=y=\sqrt{1-e^{-t}}$ and $\zeta=\cos\th$ in (\ref{eq:final Ajn}). 
\eop
Since by a simple function theoretic argument
the Legendre polynomials $P_j(\cos\th)$ on the right hand side of
(\ref{eq:Fourier by Legendre}) can be written as
\be
P_j(\cos\th)= 
\sum_{l=0}^j g_lg_{j-l}\cos (j-2l)\th
\;,
\label{eq:Pj}
\ee
with positive coefficients
\be
g_l=
\frac{(2l)!}{4^l\,l!^2}
\label{eq:glpos}
\ee
(see e.g.\ \cite{Sze39}, (4.9.3)), we have at this stage the
\bC
{\rm
The Weinstein functions satisfy the inequalities (\ref{eq:pos}),
\[
\Lambda_k^n(t)\geq 0\;\;\;\;\;\;(t\in\R^+,\;\;\;0\leq k\leq n)
\;.
\]
}
\eC
\pr
Combining Theorems~\ref{th:Weinstein Fourier}~and~\ref{th:Askey-Gasper type}
with (\ref{eq:Pj})--(\ref{eq:glpos}) gives the result.
\eop
Theorem~\ref{th:Askey-Gasper type} together with (\ref{eq:Pj}) immediately
yields sum representations for the Weinstein functions in terms of
the Gegenbauer polynomials,
\[
\Lambda_{2m}^n(t)=
\sum_{j=m}^{[n/2]}
4^{2j}\,\frac{\Gamma(n+1-2j)(2j)!^2}{\Gamma(n+2+2j)}\,(4j+1)\,g_{j-m}\,g_{j+m}\,
e^{-2jt}\lk C_{n-2j}^{2j+1}\lk\sqrt{1\!-\!e^{-t}}\rk\rk^2
\]
for $m=0,1,\ldots,[n/2]$, and
\[
\Lambda_{2m+1}^n(t)=
\!\sum_{j=m}^{[(n-1)/2]}
\!4^{2j+1}\,\frac{\Gamma(n\!-\!2j)(2j\!+\!1)!^2}{\Gamma(n+3+2j)}\,(4j+3)\,
g_{j-m}\,g_{j+1+m}\,
e^{-(2j+1)t}\lk C_{n-2j-1}^{2j+2}\lk\sqrt{1\!-\!e^{-t}}\rk\rk^2
\]
for $m=0,1,\ldots,[(n\!-\!1)/2]$. 
Another form of this statement will be given in
\S~\ref{sec:Askey-Gasper Identity for the Weinstein Functions}.

\section{Addition Theorem for the Gegenbauer Polynomials}
\label{sec:Addition Theorem for the Gegenbauer Polynomials}

In this section, we fill the gap that remained open in the
last section by proving a generalization of 
(\ref{eq:Cn1sum})/(\ref{eq:fourth step}),
the {\sl addition theorem} for the Gegenbauer polynomials (see e.g.\
\cite{Gegenbauer}).
\bT
[Addition Theorem for the Gegenbauer Polynomials]
\label{th:Addition Theorem for the Gegenbauer Polynomials}
{\rm
For $\nu>1/2$, $x,y\in [-1,1]$, and $\zeta\in \C$,
the Gegenbauer polynomials satisfy the identity
\[
C_n^\nu\Big(xy+\sqrt{1-x^2}\sqrt{1-y^2}\,\zeta\Big)
=
\]
\[
\Gamma(2\nu-1)\sum_{j=0}^n 
\frac{4^j\,(n-j)!}{\Gamma(n+2\nu+j)}\,\Big((\nu)_j\Big)^2\,(2\nu+2j-1)
(1-x^2)^{j/2}\,(1-y^2)^{j/2}\,C_{n-j}^{\nu+j}(x)\,C_{n-j}^{\nu+j}(y)\,
C_j^{\nu-1/2}(\zeta)
\;.
\]
}
\eT
\pr
The function 
\[
C_n^\nu(xy+\sqrt{1-x^2}\sqrt{1-y^2}\,\zeta\Big)
\]
as a function of $\zeta$ is a polynomial of degree $n$. Therefore, for
any $\lambda>0$, we can
expand it in terms of Gegenbauer polynomials $C_j^\lambda(\zeta)$,
\be
C_n^\nu(xy+\sqrt{1-x^2}\sqrt{1-y^2}\,\zeta\Big)
=
\sum_{m=0}^n
A_m^n(x,y)\,C_m^\lambda(\zeta)
\;,
\label{eq:Ajnansatz}
\ee
the coefficients $A_j^n$ being functions of the parameters $x$ and $y$.

The orthogonality relation of the system $C_j^\lambda(\zeta)$ is
given by
\[
\int\limits_{-1}^1 (1-\zeta^2)^{\lambda-1/2}\,
C_j^\lambda(\zeta)\,C_m^\lambda(\zeta)\,d\zeta=
\funkdef{\frac{\pi\,2^{1-2\lambda}\,\Gamma(j+2\lambda)}
{j!\,(j+\lambda)\,\Gamma(\lambda)^2}}{j=m}
{0}
\]
(see e.g.\ \cite{AS}, (22.2.3)).
Multiplying (\ref{eq:Ajnansatz})
by $(1-\zeta^2)^{\lambda-1/2}\,C_j^\lambda(\zeta)$, and
integrating from $\zeta=-1$ to $\zeta=1$, we get therefore
\[
\int\limits_{-1}^1 (1-\zeta^2)^{\lambda-1/2}\,
C_n^\nu(xy+\sqrt{1-x^2}\sqrt{1-y^2}\,\zeta\Big)
C_j^\lambda(\zeta)\,d\zeta=
A_j^n(x,y)\,
\frac{\pi\,2^{1-2\lambda}\,\Gamma(j+2\lambda)}
{j!\,(j+\lambda)\,\Gamma(\lambda)^2}
\;.
\]
Utilizing identity (\ref{eq:identity f}) with 
\[
f(\zeta):=C_n^{\nu}\Big(xy\!+\!\sqrt{1\!-\!x^2}\sqrt{1\!-\!y^2}\,\zeta\Big)
\;,
\]
we get
\[
A_j^n(x,y)=
\frac{2^{j+2\lambda-1}\,\Gamma(\lambda)\,\Gamma(j\!+\!\lambda+1)}
{\pi\,\Gamma(2j+2\lambda)}
\int\limits_{-1}^1\!
(1-\zeta^2)^{j+\lambda-1/2}\,\frac{d^j}{d\zeta^j}
C_n^{\nu}\Big(xy\!+\!\sqrt{1\!-\!x^2}\sqrt{1\!-\!y^2}\zeta\Big)\,d\zeta
\;.
\]
The derivative identity (\ref{eq:derivativeidentity}) then yields
\beao
A_j^n(x,y)&=&
\frac{2^{2j+2\lambda-1}\,(\nu)_j\,\Gamma(\lambda)\,\Gamma(j+\lambda+1)}
{\pi\,\Gamma(2j+2\lambda)}
(1-x^2)^{j/2}\,(1-y^2)^{j/2}
\\&&\cdot
\int\limits_{-1}^1
(1-\zeta^2)^{j+\lambda-1/2}\,
C_{n-j}^{\nu+j}\Big(xy+\sqrt{1-x^2}\sqrt{1-y^2}\zeta\Big)\,d\zeta
\;.
\eeao
Now we choose $\lambda:=\nu-1/2$ (hence our assumption $\nu>1/2$).
This choice is motivated by the
calculation involving the differential equation that follows later, for
which the desired simplification occurs exactly when $\lambda=\nu-1/2$.
Using the duplication formula 
\be
\Gamma(2z)=\frac{2^{2z-1}}{\sqrt{\pi}}\,\Gamma(z)\,\Gamma(z+1/2)
\label{eq:duplication formula}
\ee
of the Gamma function to simplify the factor in front of the integral, we
finally arrive at the representation 
\beao
A_j^n(x,y)&=&
2^{1-2\nu}(2j+2\nu-1)\frac{\Gamma(2\nu-1)}{\Gamma(\nu)^2}\,
(1-x^2)^{j/2}\,(1-y^2)^{j/2}
\\&&\cdot
\int\limits_{-1}^1
(1-\zeta^2)^{j+\nu-1}\,
C_{n-j}^{\nu+j}\Big(xy+\sqrt{1-x^2}\sqrt{1-y^2}\zeta\Big)\,d\zeta
\eeao
for the coefficients $A_j^n(x,y)$.
Hence, we consider the function
\[
Q_j^n(x,y):=
\int\limits_{-1}^1
(1-\zeta^2)^{j+\nu-1}\,
C_{n-j}^{\nu+j}\Big(xy+\sqrt{1-x^2}\sqrt{1-y^2}\zeta\Big)\,d\zeta
\]
in detail.
Observe that $Q_j^n(x,y)$ is a polynomial in the variables $x$ and $y$,
of degree $n-j$ each. Note furthermore that $Q_j^n(x,y)$ is symmetric,
i.e.\ $Q_j^n(x,y)=Q_j^n(y,x)$. 

In the following we will show that the
integral $Q_j^n(x,y)$ has zeros at both the zeros
of $C_{n-j}^{\nu+j}(x)$ and $C_{n-j}^{\nu+j}(y)$, hence, as a polynomial
of degree $n-j$ in $x$ and $y$ respectively, must be a 
constant multiple of the product
$C_{n-j}^{\nu+j}(x)\,C_{n-j}^{\nu+j}(y)$. 

By the symmetry of $Q_j^n(x,y)$
it is enough to show that $Q_j^n(x,y)$ has zeros at the zeros
of $C_{n-j}^{\nu+j}(x)$. Since $C_{n-j}^{\nu+j}(x)$ is a solution of the
differential equation
\be
(1-x^2)\,p''(x)-(2\nu+2j+1)\,x\,p'(x)+(n-j)(n+j+2\nu)\,p(x)=0
\;,
\label{eq:DE}
\ee
and since any {\sl polynomial} solution $p(x)$ of (\ref{eq:DE}) must be a 
multiple of $C_{n-j}^{\nu+j}(x)$ (see e.g.\ \cite{Sze39}, Theorem 4.2.2
in combination with \cite{AS}, (22.5.27)),
we have only to check that $p(x):=Q_j^n(x,y)$ satisfies (\ref{eq:DE}).

We write $\eta(x):=xy+\sqrt{1-x^2}\sqrt{1-y^2}\zeta$, and note that
\[
\eta'(x)=y-\frac{\sqrt{1-y^2}}{\sqrt{1-x^2}}\,x\,\zeta
\]
so that
\[
x\eta'(x)=xy-\frac{\sqrt{1-y^2}}{\sqrt{1-x^2}}\,x^2\,\zeta
=\eta(x)-\frac{\sqrt{1-y^2}}{\sqrt{1-x^2}}\,\zeta
\;.
\]
Hence we deduce
\beao
-(2\nu+2j+1)\,x\,\ded x Q_j^n(x,y)
&=&
\int\limits_{-1}^1
-(2\nu+2j+1)\,\eta(x)\,
\Big( C_{n-j}^{\nu+j}\Big)'(\eta(x))
(1-\zeta^2)^{j+\nu-1}\,
d\zeta
\\&&
+
(2\nu+2j+1)\frac{\sqrt{1-y^2}}{\sqrt{1-x^2}}
\int\limits_{-1}^1
\zeta\,(1-\zeta^2)^{j+\nu-1}\,
\Big( C_{n-j}^{\nu+j}\Big)'(\eta(x))\,d\zeta
.
\eeao
Similarly, using the identity
\[
\lk y\,\sqrt{1-x^2}-x\,\sqrt{1-y^2}\,\zeta\rk^2=
(1-\eta(x)^2)-(1-y^2)\,(1-\zeta^2)
\;,
\]
we get
\beao
(1-x^2)\,\dedn{x}{2} Q_j^n(x,y)
&=&
\int\limits_{-1}^1
(1-\eta(x)^2)\,
\Big( C_{n-j}^{\nu+j}\Big)''(\eta(x))(1-\zeta^2)^{j+\nu-1}\,
d\zeta
\\&&
-
\frac{\sqrt{1-y^2}}{\sqrt{1-x^2}}
\int\limits_{-1}^1
\sqrt{1-x^2}\sqrt{1-y^2}\,(1-\zeta^2)^{j+\nu}\,
\Big( C_{n-j}^{\nu+j}\Big)''(\eta(x))\,d\zeta
\\&&
-
\frac{\sqrt{1-y^2}}{\sqrt{1-x^2}}
\int\limits_{-1}^1
\zeta\,(1-\zeta^2)^{j+\nu-1}\,
\Big( C_{n-j}^{\nu+j}\Big)'(\eta(x))\,d\zeta
.
\eeao
Combining these results, we arrive at the representation
\[
(1-x^2)\,\dedn{x}{2} Q_j^n(x,y)-(2\nu+2j+1)\,x\,\ded{x} Q_j^n(x,y)+
(n-j)(n+j+2\nu)\,Q_j^n(x,y)=
\]
\[
\int\limits_{-1}^1
(1\!-\!\zeta^2)^{j+\nu-1}\,
\!\!\left(
(1\!-\!\eta^2)\,\Big( C_{n-j}^{\nu+j}\Big)''(\eta)
\!-\!(2\nu\!+\!2j\!+\!1)\,\eta
\Big( C_{n-j}^{\nu+j}\Big)'(\eta)\!+
\!(n\!-\!j)(n\!+\!j\!+\!2\nu)\,C_{n-j}^{\nu+j}(\eta)
\right)
\!d\zeta
\]
\[
+
\frac{\sqrt{1\!-\!y^2}}{\sqrt{1\!-\!x^2}}
\lk
\int\limits_{-1}^1
2(j\!+\!\nu)\zeta(1\!-\!\zeta^2)^{j+\nu-1}
\Big( C_{n-j}^{\nu+j}\Big)'(\eta)d\zeta
\!-
\!\int\limits_{-1}^1
\!\!\sqrt{1\!-\!x^2}\sqrt{1\!-\!y^2}(1\!-\!\zeta^2)^{j+\nu}
\Big( C_{n-j}^{\nu+j}\Big)''(\eta)d\zeta
\rk
\!\!.
\]
The first integral obviously vanishes since $C_{n-j}^{\nu+j}(x)$ satisfies
the differential equation
(\ref{eq:DE}). The vanishing of the final parenthesized expression 
follows easily by partial integration. Therefore, we have proved that
$Q_j^n(x,y)$ is a solution of (\ref{eq:DE}), as announced.

Hence, 
\be
Q_j^n(x,y)=a\,C_{n-j}^{\nu+j}(x)\,C_{n-j}^{\nu+j}(y)
\label{eq:Bestimmung von a}
\ee
with a constant $a$ (not depending on $x$ and $y$). For $y=1$, we deduce
\be
Q_j^n(x,1)=
\int\limits_{-1}^1
(1-\zeta^2)^{j+\nu-1}\,
C_{n-j}^{\nu+j}(x)\,d\zeta
=
2^{2j+2\nu-1}\frac{\Gamma(j+\nu)^2}{\Gamma(2j+2\nu)}\,
C_{n-j}^{\nu+j}(x)
\label{eq:Beta evaluation}
\ee
by an evaluation of the Beta type integral. On the other hand, by
(\ref{eq:Bestimmung von a}),
\[
Q_j^n(x,1)=a\,C_{n-j}^{\nu+j}(x)\,C_{n-j}^{\nu+j}(1)
=
a\,C_{n-j}^{\nu+j}(x)\,\ueber{n+j+2\nu-1}{n-j}
\]
(see e.g.\ \cite{AS}, (22.4.2)), so that we get 
\[
a=2^{2j+2\nu-1}\frac{\Gamma(j+\nu)^2}{\Gamma(2j+2\nu)}\left/
\ueber{n+j+2\nu-1}{n-j}\right.
=
2^{2j+2\nu-1}\frac{(n-j)!\,\Gamma(j+\nu)^2}{\Gamma(n+j+2\nu)}
\;.
\]
Hence
\[
Q_j^n(x,y)=
2^{2j+2\nu-1}\frac{(n-j)!\,\Gamma(j+\nu)^2}{\Gamma(n+j+2\nu)}
\,
C_{n-j}^{\nu+j}(x)\,C_{n-j}^{\nu+j}(y)
\;,
\]
implying
\[
A_j^n(x,y)=
\Gamma(2\nu-1)\frac{2^{2j}(n-j)!}{\Gamma(n\!+\!j\!+\!2\nu)}\,
\frac{\Gamma(j\!+\!\nu)^2}{\Gamma(\nu)^2}\,
(2j+2\nu-1)
\,
(1-x^2)^{j/2}\,(1-y^2)^{j/2}
\,
C_{n-j}^{\nu+j}(x)\,C_{n-j}^{\nu+j}(y)
,
\]
and we are done.
\eop
As a consequence, taking the limit $\nu\pf 1/2$, we get the following
\bC
[Addition Theorem for the Legendre Polynomials]
{\rm
For $x,y\in [-1,1]$, $\zeta\in \C$,
the Legendre polynomials satisfy the identities
\[
P_n(xy+\sqrt{1-x^2}\sqrt{1-y^2}\,\zeta\Big)
=
\]
\be
P_n(x)\,P_n(y)+
2\sum_{j=1}^n 4^j\,\frac{(n-j)!}{(n+j)!}\Big( (1/2)_j\Big)^2\,
(1-x^2)^{j/2}\,(1-y^2)^{j/2}\,C_{n-j}^{1/2+j}(x)\,C_{n-j}^{1/2+j}(y)\,
T_j(\zeta)
\label{eq:CjTj}
\ee
\be
=
P_n(x)\,P_n(y)+
2\sum_{j=1}^n \frac{(n-j)!}{(n+j)!}\,
P_n^j(x)\,P_n^j(y)\,T_j(\zeta)
\;,
\label{eq:PjTj}
\ee
where $T_j(\zeta)$ denote the Chebyshev polynomials of the first kind, and
\be
P_n^j(x)=(-1)^j\,(1-x^2)^{j/2}\,\dedn{x}{j} P_n(x)
\label{eq:Pnj}
\ee
denote the associated Legendre functions (see e.g.\ \cite{AS}, (8.6.6)).

In particular, for $y=x$, one has
\be
P_n(x^2+(1-x^2)\cos\th )=
P_n(x)^2+
2\sum_{j=1}^n \frac{(n-j)!}{(n+j)!}\,
P_n^j(x)^2\,\cos j\th
\;.
\label{eq:PnWeinstein}
\ee
}
\eC
\pr
Since 
\[
C_n^0(x)=\lim_{\lambda\pf 0}\frac{C_n^\lambda(x)}{\lambda}
\quad\quad\mbox{and}\quad\quad
C_n^\alpha(x)=\lim_{\lambda\pf\alpha} C_n^\lambda(x)
\quad
\mbox{for all }\alpha>0
\]
(see e.g.\ \cite{AS}, (22.5.4)), for $\nu\pf1/2$
Theorem~\ref{th:Addition Theorem for the Gegenbauer Polynomials} implies
\[
C_n^{1/2}\Big(xy+\sqrt{1-x^2}\sqrt{1-y^2}\,\zeta\Big)
\]
\[
=
C_n^{1/2}(x)\,C_n^{1/2}(y)
+\sum_{j=1}^n 
4^j\,\frac{(n-j)!}{(n+j)!}\,\Big((1/2)_j\Big)^2\,
(1-x^2)^{j/2}\,(1-y^2)^{j/2}\,C_{n-j}^{1/2+j}(x)\,C_{n-j}^{1/2+j}(y)\,
j\,C_j^{0}(\zeta)
\;.
\]
With $C_n^{1/2}(x)=P_n(x)$, and $j\,C_j^{0}(\zeta)=2\,T_j(\zeta)$
(see e.g.\ \cite{AS}, (22.5.35), (22.5.33)), we
get (\ref{eq:CjTj}). An application of (\ref{eq:derivativeidentity})
and (\ref{eq:Pnj}) yields (\ref{eq:PjTj}).

Using 
\[
T_n(\cos\th)=\cos n\th
\]
(see e.g.\ \cite{AS}, (22.3.15)) finally yields (\ref{eq:PnWeinstein}).
\eop
Note that Weinstein used (\ref{eq:PnWeinstein}) in his proof of Milin's
conjecture.

\section{Askey-Gasper Identity for the Weinstein Functions}
\label{sec:Askey-Gasper Identity for the Weinstein Functions}

Here, we combine the above results to deduce a sum representation
with nonnegative summands for the Weinstein functions, and therefore by 
Theorem~\ref{th:Weinstein Jacobisum} for the Jacobi polynomial sum.

By Theorem~\ref{th:Askey-Gasper type} we have
\[
C_n^{1}\Big((1-e^{-t})+e^{-t}\cos\th\Big)
=
\sum_{j=0}^n
(2j+1)
\frac{4^{j}\,j!^2\,(n-j)!}{(n+j+1)!}\,
e^{-jt}\,\lk C_{n-j}^{j+1}\left(\sqrt{1-e^{-t}}\right)\rk^2\,
P_j(\cos\th)
\;,
\]
and, expanding $P_j(\cos\th)$ using (\ref{eq:CjTj}) with $x=y=0$, 
$\zeta=\cos\th$, this gives
\[
=
\sum_{j=0}^n
(2j+1)
\frac{4^{j}j!^2(n\!-\!j)!}{(n+j+1)!}
e^{-jt}\lk C_{n-j}^{j+1}\left(\sqrt{1\!-\!e^{-t}}\right)\rk^2\cdot
2\sumprime_{k=0}^j 4^k\frac{(j\!-\!k)!}{(j\!+\!k)!}\Big( (1/2)_k\Big)^2
C_{j-k}^{1/2+k}(0)^2
T_k(\cos\th)
,
\]
where $\Sigma'$ indicates that the summand for $k=0$ is to be taken with a
factor $1/2$.
Interchanging the order of summation, and using $T_k(\cos\th)=\cos k\th$, gives
\[
=
2\sumprime_{k=0}^n\sum_{j=k}^n
(2j+1)4^k
\frac{4^{j}j!^2(n\!-\!j)!}{(n+j+1)!}
\frac{(j\!-\!k)!}{(j\!+\!k)!}\Big( (1/2)_k\Big)^2
e^{-jt}C_{j-k}^{1/2+k}(0)^2
\lk C_{n-j}^{j+1}\left(\sqrt{1\!-\!e^{-t}}\right)\rk^2
\cos k\th
.
\]
Comparing with Theorem~\ref{th:Weinstein Fourier}, 
\[
C_n^{1}\Big((1-e^{-t})+e^{-t}\cos\th\Big)
=
2\sumprime_{k=0}^n \Lambda_k^n(t) \cos k\th
\;,
\]
and equating coefficients yields for the Weinstein functions
\[
\Lambda_k^n(t)=
\sum_{j=k}^n
(2j+1)4^k
\frac{4^{j}j!^2(n\!-\!j)!}{(n+j+1)!}
\frac{(j\!-\!k)!}{(j\!+\!k)!}\Big( (1/2)_k\Big)^2
e^{-jt}C_{j-k}^{1/2+k}(0)^2
\lk C_{n-j}^{j+1}\left(\sqrt{1\!-\!e^{-t}}\right)\rk^2
\;.
\]
Replacing $n$ by $k+n$, and then making the index shift 
$j_{\rm new}:=j_{\rm old}-k$ finally leads to
\[
\Lambda_k^{k+n}(t)\!=
\!\sum_{j=0}^{n}
(2j\!+\!2k\!+\!1)
\frac{4^{j+2k}(j\!+\!k)!^2(n\!-\!j)!j!( (1/2)_k)^2}
{(2k+n+j+1)!(j+2k)!}
e^{-(j+k)t}C_{j}^{1/2+k}(0)^2
{\lk C_{n-j}^{j+k+1}\!\!\left(\sqrt{1\!-\!e^{-t}}\right)\rk\!}^2
\!.
\]
Setting $y:=\sqrt{1\!-\!e^{-t}}$, by Theorem~\ref{th:Weinstein Jacobisum}
\[
\sum_{j=0}^{n} \!P_j^{(2k,0)}(2y^2\!-\!1)\!=
\!\sum_{j=0}^{n}
(2j\!+\!2k\!+\!1)
\frac{4^{j+2k}(j\!+\!k)!^2(n\!-\!j)!j!( (1/2)_k)^2}
{(2k+n+j+1)!(j+2k)!}
(1\!-\!y^2)^jC_{j}^{1/2+k}(0)^2
{\lk C_{n-j}^{j+k+1}(y)\!\rk\!}^2
\!.
\]
This is an Askey-Gasper type representation different from
(\ref{eq:AG Identity}) that was given by Gasper
(\cite{Gasper}, (8.17), and (8.18) with $x=0$). 
Note that Gasper's formula 
(\cite{Gasper}, (8.18)) interpolates between these two representations.
Whereas Askey's and Gasper's 
deductions of the given formulas prove the results for all $\alpha>-2$,
our deduction has the disadvantage that it is only valid for 
$\alpha=2k, k\in\N_0$.
On the other hand, the advantage of our presentation is that it embeds
this result in a natural way in Weinstein's proof of Milin's conjecture
using only elementary properties of classical orthogonal polynomials.

\section{Closed Form Representation of \mbox{Weinstein functions}}
\label{sec:Closed Form Representation of Weinstein functions}

Note that nowhere in our deduction we needed the explicit representation
of the de Branges functions = Weinstein functions, compare 
Henrici's comment \cite{Henrici}, p.\ 602: ``At the time of this writing,
the only way to verify 
$\dot{\tau_{k}^n}(t)\leq 0$ appears to be to
solve the system explicitly, and to manipulate the solution''.

In this connection we like to mention
that in \cite{KSdeBranges} we proved the identity
(\ref{eq:deBrangesWeinstein}), which connects de Branges' with Weinstein's
functions, by a pure application of the de Branges differential
equations system (see also \cite{Schmersau}), 
and without the use of an explicit representation of the de Branges functions.

In this section we give
a simple method to generate this explicit
representation which was used by de Branges, see also \cite{Wilf}.

Since $(1-e^{-t})+e^{-t}\cos\th=1-2e^{-t}\sin^2\frac{\th}{2}$,
Taylor expansion gives using (\ref{eq:derivativeidentity}) and
(\cite{AS}, (22.4.2))
\beao
C_n^{1}\Big((1-e^{-t})+e^{-t}\cos\th\Big)
&=&
C_n^{1}\Big(1-2e^{-t}\sin^2\frac{\th}{2}\Big)
=
\sum_{j=0}^n 
\frac{{C_n^{1}}^{(j)}(1)}{j!}\,(-1)^j\,2^j\,e^{-jt}
\lk\sin^{2}\frac{\th}{2}\rk^j
\\&=&
\sum_{j=0}^n 
C_{n-j}^{j+1}(1)\,2^{2j}\,(-1)^j\,e^{-jt}
\lk\sin^{2}\frac{\th}{2}\rk^j
\\&=&
\sum_{j=0}^n
\ueber{n+j+1}{n-j}\,2^{2j}\,(-1)^j\,e^{-jt}
\lk\sin^{2}\frac{\th}{2}\rk^j
\;.
\eeao
An elementary argument shows that
\[
\lk\sin^{2}\frac{\th}{2}\rk^j
=
2\sumprime_{k=0}^j 
\frac{(-1)^k}{2^{2j}}\,\ueber{2j}{j-k} T_{2k}\lk\cos \frac{\th}{2}\rk
=
2\sumprime_{k=0}^j \frac{(-1)^k}{2^{2j}}\,\ueber{2j}{j-k}\cos k\th
\]
(see e.g.\ \cite{Tri}, p.\ 189). Changing the order of summation, we get
therefore
\beao
C_n^{1}\Big((1-e^{-t})+e^{-t}\cos\th\Big)
&=&
2\sumprime_{k=0}^n\sum_{j=k}^n
(-1)^{j+k}\,\ueber{n+j+1}{n-j}\,\ueber{2j}{j-k}\,e^{-jt}\,\cos k\th
\\&=&
2\sumprime_{k=0}^n \Lambda_k^n(t) \cos k\th
\eeao
by (\ref{eq:Weinstein Fourier}). Hence
\beao
\Lambda_k^n(t)
&=&
\sum_{j=k}^n (-1)^{j+k}\,\ueber{n+j+1}{n-j}\,\ueber{2j}{j-k}\,e^{-jt}
\\&=&
e^{-kt}\,\ueber{n+k+1}{n-k}\;
\hypergeom{3}{2}{n+k+2,k+1/2,-n+k}{k+3/2,2k+1}{e^{-t}}
\;.
\eeao

\section*{Acknowledgement}
The first author
would like to thank Peter Deuflhard who initiated his studies on the 
work with orthogonal polynomials.


\begin{thebibliography}{99}

\bi{AS}
Abramowitz, M.\ and Stegun, I.\ A.: {\sl Handbook of Mathematical
Functions.} Dover Publ., New York, 1964.

\bi{AG}
Askey, R. and Gasper, G.: Positive Jacobi polynomial sums II. Amer.\
J.\ Math.\ \9\8 (1976), 709--737.


\bi{Bieberbach}
Bieberbach, L.: \"Uber die Koeffizienten derjenigen Potenzreihen,
welche eine schlichte Abbildung des Einheitskreises vermitteln. S.-B.\
Preuss.\ Akad.\ Wiss.\ {\bf 38} (1916), 940--955.

\bi{Bra}
De Branges, L.: A proof of the Bieberbach conjecture. Acta Math. \1\5\4
(1985), 137--152.



\bi{FP}
FitzGerald, C.\ H. and Pommerenke, Ch.: The de Branges Theorem on
univalent functions. Trans.\ Amer.\ Math.\ Soc.\ \2\9\0 (1985), 683--690.

\bi{Gasper}
Gasper, G.: Positivity and special functions. In: {\sl Theory and Application
of Special Functions.} Edited by R.A.\ Askey. Academic Press, New York, 1975,
375--433.

\bi{Gegenbauer}
Gegenbauer, L.: Das Additionstheorem der Funktionen $C_n^\nu(x)$.
Sitzungsberichte der mathematisch-naturwissenschaftlichen Klasse
der Akademie der Wissenschaften Wien Abteilung II a, \1\0\2 (1893),
942--950.


\bi{Henrici}
Henrici, P.: {\sl Applied and Computational Complex Analysis, Vol.\ 3:\
Discrete Fourier Analysis -- Cauchy Integrals -- Construction of
Conformal maps -- Univalent Functions}. John Wiley \& Sons, New York,
1986.


\bi{Hua}
Hua, L.K.: {\sl Harmonic Analysis of Functions of Several Complex Variables
in the Classical Domains.} Translations of Mathematical Monographs Vol.\ 6,
Amer.\ Math.\ Soc., Providence, R.I., 1963.

\bi{Koepf}
Koepf, W.:
Von der Bieberbachschen Vermutung zum Satz von de Branges sowie
der Beweisvariante von Weinstein. In: 
{\sl Jahrbuch \"Uberblicke Mathematik 1994}. 
Vieweg-Verlag, Braunschweig--Wiesbaden, 1994, 175--193.



\bi{KSdeBranges}
Koepf, W. and Schmersau, D.:
On the de Branges theorem.
Konrad-Zuse-Zentrum Berlin (ZIB), Preprint SC 95-10, 1995.

\bi{LM}
Lebedev, N.\ A.\ and Milin, I.\ M.: An inequality. Vestnik Leningrad
Univ.\  \2\0 (1965), 157--158 (Russian).


\bi{Loewner2}
L\"owner, K.: Untersuchungen \"uber schlichte konforme Abbildungen
des Einheitskreises I. Math.\ Ann.\ \8\9 (1923), 103--121.


\bi{Schmersau}
Schmersau, D.: Untersuchungen zur Rekursion von L.\ de Branges.
Complex Variables \1\5 (1990), 115--124.

\bi{Sze39}
Szeg\"o, G.: {\sl Orthogonal Polynomials}. Amer.\ Math.\ Soc.\ Coll.\
Publ.\ Vol.\ \2\3, New York City, 1939.

\bi{Todorov}
Todorov, P.: A simple proof of the Bieberbach conjecture.
Bull.\ Cl.\ Sci., VI.\ S\'er., Acad.\ R.\ Belg.\ 3 \1\2 (1992), 335--356.

\bi{Tri}
Tricomi, F.\ G.: {\sl Vorlesungen \"uber Orthogonalreihen}. Grundlehren
der Mathematischen Wissenschaften \7\6, Springer-Verlag,
Berlin--G\"ottingen--Heidelberg, 1955.

\bi{Weinstein}
Weinstein, L.: The Bieberbach conjecture.
International Mathematics Research Notices \5 (1991), 61--64.

\bi{Wilf}
Wilf, H.: A footnote on two proofs of the Bieberbach-de Branges Theorem.
Bull.\ London Math.\ Soc.\ \2\6 (1994), 61--63.


%
\end{thebibliography}
\end{document}